\documentclass[10pt]{article}
\usepackage{amsmath,amssymb,amsthm,amscd}
\numberwithin{equation}{section}

\def\cD{\mathcal D}
\def\cE{{\mathcal E}}

\def\cH{{\mathcal H}}

\def\cM{{\mathcal M}}

\def\cR{{\mathcal R}}

\def\cU{{\mathcal U}}
\def\cV{{\mathcal V}}
\def\cW{{\mathcal W}}

\def\bE{{\mathbf E}}

\def\RR{{\mathbb R}}
\def\CC{{\mathbb C}}

\newtheorem{prop}{Proposition}[section]
\newtheorem{theo}[prop]{Theorem}
\newtheorem{lemm}[prop]{Lemma}
\newtheorem{coro}[prop]{Corollary}

\newtheorem{defi}[prop]{Definition}

\def\begeq{\begin{equation}}
\def\endeq{\end{equation}}

\def\and{\quad{\rm and}\quad}

\def\lab{\ }
\def\lab{\label}

\title{Geometry of K\"ahler metrics and holomorphic foliation by discs}

\author{X. X. Chen and G. Tian\thanks{Both authors are supported by
NSF research grants and the second author is also supported
partially by a Simons fund}}
\date{}

\begin{document}
\bibliographystyle{plain}

\maketitle
\section{Introduction and Main Results}

The purpose of this paper is to establish a partial regularity
theory on certain homogeneous complex Monge-Ampere equations. As
consequences of this new theory, we prove the uniqueness of
extremal K\"ahler metrics and give an necessary condition for
existence of extremal K\"ahler metrics.

Following \cite{calabi82}, we call a K\"ahler metric extremal if
the complex gradient of its scalar curvature is a holomorphic
vector field. In particular, any K\"ahler metric with constant
scalar curvature is extremal, conversely, if the underlying
K\"ahler manifold has no holomorphic vector fields, then an
extremal K\"ahler metric is of constant scalar curvature. Our
first result is
\begin{theo}
\lab{th:unique} Let $(M, [\omega])$ be a compact K\"ahler manifold
with a K\"ahler class $[\omega]\in H^2(M,\RR)\cap H^{1,1}(M,\CC)$.
Then there is at most one extremal K\"ahler metric with K\"ahler
class $[\omega]$ modulo holomorphic transformations, that is, if
$\omega_1$ and $\omega_2$ are extremal K\"ahler metrics with the
same K\"ahler class, then there is a holomorphic transformation
$\sigma$ such that $\sigma^*\omega_1 =\omega_2$.
\end{theo}

The problem of uniqueness of extremal K\"ahler metrics has a long
history. The uniqueness of K\"ahler-Einstein metrics was pointed
out by Calabi in early 50's in the case of non-positive scalar
curvature. In \cite{Bando87}, Bando and Mabuchi proved that the
uniqueness of K\"ahler-Einstein metric in the case of positive
scalar curvature.
% In \cite{tianzhu00} and \cite{tianzhu01}, G.Tian and X. H. Zhu proved the uniqueness of K\"{a}hler-Ricci Solitons.
 Following a suggestion of Donaldson, the first author
proved in \cite{chen991} the uniqueness of K\"ahler metrics with
constant scalar curvature in any K\"ahler class which admits a
K\"ahler metric with non-positive scalar curvature. In
\cite{dona02}, S. Donaldson proved the uniqueness of constant
scalar curvature K\"ahler metrics with rational K\"ahler class on
any projective manifolds (which are K\"ahler) without non-trivial
holomorphic vector fields.\footnote{After we finished proving our
uniqueness theorem, we learned that T. Mabuchi extended S.
Donaldson's arguements to any extremal K\"ahler metrics with
rational coefficients on any projective manifolds and proved their
uniqueness in those special cases.}

The existence of extremal metrics remains open in general cases.
One difficulty is due to the fact the associated equation is fully
non-linear and of 4th order. The case of K\"ahler-Einstein metrics
has been well understood (see \cite{Yau78}, \cite{Au76},
\cite{Tian97}).

As another consequence of our regularity theory, we will give an
necessary condition on existence of K\"ahler metrics with constant
scalar curvature in terms of the K-energy. In \cite{Ma87}, Mabuchi
introduced the following functional $\bE_\omega$ as follows: For
any $\varphi$ with $\omega_\varphi =
\omega+\partial\overline{\partial}\varphi > 0$, define $$
\bE_\omega(\varphi)=-\int_0^1\int_M
\dot\varphi(s(\omega_{\varphi_t})-\mu)\omega_{\varphi_t}^n\wedge
dt,$$ where $\omega_{\varphi_t}$ is any path of K\"ahler metrics
joining $\omega$ and $\omega_\varphi$, $s(\omega_{\varphi_t})$
denotes the scalar curvature and $\mu $ is its average. It turns
out that $\mu$ is determined by the first Chern class $c_1(M)$ and
the K\"ahler class $[\omega]$.

\begin{theo}
\label{th:minimum} Let $M$ be a compact K\"ahler manifold with a
constant scalar curvature K\"ahler metric $\omega$. Then
$\bE_\omega(\varphi)\ge 0$ for any $\varphi$ with $\omega_\varphi
>0$.
\end{theo}

This theorem was proved for K\"ahler-Einstein metrics in
\cite{Bando87} (also see \cite{Tian97}) and in \cite{chen991} for
K\"ahler manifolds with non-positive first Chern class. This
theorem can be also generalized to arbitrary extremal K\"ahler
metrics by using the modified K-energy. This theorem gives a
partial answer to a conjecture of the second author: $M$ has a
constant scalar curvature K\"ahler metric in a given K\"ahler
class $[\omega]$ if and only if the K-energy is proper in a
suitable sense on the space of K\"ahler metrics with the fixed
K\"ahler class $[\omega]$. We will further discuss applications of
our method here to this problem on properness in a forthcoming
paper. Combining Theorem 1.2 with results in \cite{tian001} and
\cite{paultian04}, we can deduce

\begin{coro} Let $(M,L)$ be a polarized algebraic manifold, that is,
$M$ is algebraic and $L$ is a positive line bundle. If there is a
constant scalar curvature K\"ahler metric with K\"ahler class
equal to $c_1(L)$. Then $(M,L)$ is asymptotically K-semistable or
CM-semistable in the sense of \cite{Tian97} (also see
\cite{tian001}) \footnote{According to \cite{paultian04}, the
CM-stability (semistability) is equivalent to the K-stability
(semistability).}.
\end{coro}

The proof of these two theorems is based on studying certain
homogenous Complex Monge-Ampere equations. Deep works have been
done on these equations (cf. \cite{Semmes92}, \cite{Dona96},
\cite{chen991}). They are related to the geodesic equation on the
space of K\"ahler metrics with $L^2$-metric.

Let $M$ be a compact K\"ahler manifold with a fixed K\"ahler
metric $\omega$ and $\Sigma$ be a Riemann surface with boundary
$\partial \Sigma$. Consider
\begin{equation} \lab{eq:hcma1}
(\pi_2^*\omega +\partial\overline{\partial }\phi)^{n+1} = 0~~{\rm
on}~\Sigma\times M,~~\phi|_{\partial \Sigma \times M}= \psi,
\end{equation}
where $\pi_2:\Sigma\times M\to M$ is the projection and $\phi$ is
a function on $\Sigma\times M$ such that $\phi (z,\cdot)\in
\cH_\omega$ for any $ z\in \Sigma$, and $\psi$ is a given function
on $\partial \Sigma\times M$ such that $\psi(z,\cdot\in \cH$. Here
$\cH_\omega$ denotes the space of K\"ahler potentials
\begin{equation}
\lab{eq:kahlerpotentials} \cH_\omega = \{ \varphi\in
C^\infty(M,\RR)~\mid ~\omega_{\varphi} = \omega +
\partial\overline{\partial} \varphi > 0,\;{\rm on} \; M\}.
\end{equation}

The proof of Theorem 1.1 and 1.2 starts with the following
observations: Given two functions $\varphi_0$ and $\varphi_1$ in
$\cH_\omega$, if there is a bounded smooth solution $\phi$ of
(\ref{eq:hcma1}) on $[0,1]\times\RR \times M$, then evaluation
function $f=\bE (\phi(z,\cdot))$ is a bounded subharmonic function
which is constant along each boundary component of $[0,1]\times
\RR,$ \footnote{This follows from the convexity of the K-energy
$\bE$} furthermore, if $\varphi_0$ is a critical metric of $\bE$,
then it follows from the Maximum principle that $\bE(\varphi_1)\ge
\bE(\varphi_0)$ and equality holds if and only if each
$\phi(z,\cdot)$ is a critical metric of $\bE$. The infinite strip
$[0,1]\times \RR$ can be approximated by discs
$\Sigma_R=[0,1]\times [-R, R]$ ($R\to \infty$). Hence, if we can
show that the equation (\ref{eq:hcma1}) has a uniformly bounded
solution for each $\Sigma_R$\footnote{$\Sigma_R$ has four corners,
but we can easily smooth corners to get a smooth Riemann surface
of disc type and with boundary and use smoothed ones to
approximate the given infinite strip.}, then Theorem 1.1 and 1.2
follow.

However, it is an extremely difficult problem to solve degenerate
complex Monge-Ampere equations. Higher regularity was obtained by
L. Cafferali, J. Kohn, L. Nirenberg and J. Spruck in 80's for
nondegenerate complex Monge-Ampere equations under certain
convexity assumptions. Weak solutions for homogenous complex Monge
Ampere equations, say in $L^p$ or $W^{1,p}$-norms, were
extensively studied (cf. \cite{Bedford76}). In \cite{chen991}, the
first named author proved the following theorem, which plays a
fundamental role in this paper.

\begin{theo}(\cite{chen991}) For any smooth map
$\psi: \partial \Sigma \rightarrow \cH$, (\ref{eq:hcma1}) always
has a unique $C^{1,1}$-solution $\phi$ on $\Sigma \times M $ such
that $\phi = \psi$ along $\partial \Sigma$. \footnote{It is not
known if $\phi(z,\cdot)$ lies in $\cH$, but the first author
proved that $\phi(z,\cdot)$ is always the limit of functions in
$\cH$.} Moreover, the $C^{1,1}$ bound of $\phi$ depends only on
the $C^2$ bound of $\psi$.
\end{theo}

Our new technical contribution is to establish partial regularity
of solutions from the above theorem in the case of $\Sigma$ being
a disc. Precisely, we prove

\begin{theo}
\lab{th:almostsmoothsolution} Let $\Sigma$ be a holomorphic disc.
For a generic boundary map $\psi: \partial \Sigma \rightarrow
\cH_\omega$, there exists a unique $C^{1,1}$ solution $\phi$ of
(\ref{eq:hcma1}) with the following properties: There is an open
and dense subset ${\cal R}_\phi\subset\Sigma\times M$ such that
\begin{enumerate}
\item $\cR_\phi$ is open and dense in $\Sigma\times M,$
% but dense in $\partial \Sigma \times M$
and the varying volume form $\omega_{\phi(z,\cdot)}^n$ extends to
a continuous function on $\Sigma_0\times M$, where
$\Sigma_0=(\Sigma\backslash \partial \Sigma).$  Moreover, it is
positive in ${\cal R}_\phi$ and vanishes identically on its
complement;
 %\item it is partially smooth,
%\item   each $\cS_\phi \cap (\{z\}\times M)$ has
%codimension at least one for each $z\in \Sigma,\;$
\item The distribution ${\cal D}_\phi$(cf. equation
(\ref{eq:distribution})) extends to a continuous distribution in a
open saturated\footnote{Any maximal extension of the leaf vector
field lies completely inside $\tilde \cV_\phi $.} set $\tilde
\cV\subset \Sigma \times M$, such that the complement $S_\phi$ of
$\tilde \cV_\phi$ is locally extendable\footnote{A closed subset
$S\subset \Sigma \times M $ of measure $0$ is {\it locally
extendable} if for any continuous function in $\Sigma\times M$
which is $C^{1,1}$ on $\Sigma\times M \setminus S$ can be extended
to a $C^{1,1}$ function on $\Sigma\times M$. Notice that any set
of codimension 2 or higher is automatically locally extendable.}
and $\phi$ is $C^1$ continuous on $\tilde \cV$. The set $S_\phi$
is referred as the singular set of $\phi$. \item The leaf vector
field ${\partial \over {\partial z}} + v$ of $\cD_\phi$ is
uniformly bounded in $\tilde \cV_\phi.\;$
\end{enumerate}
\end{theo}

We will call the solution in the above theorem an almost smooth
solution. The partial regularity in Theorem 1.5 is sharp since we
do have examples where the solution for (\ref{eq:hcma1}) is
singular. It seems to be the first time to use singular foliations
systematically to study partial regularity for homogeneous complex
Monge-Ampere equations.  Now let us explain briefly how this
theorem is proved.

It has been known for long that solutions of homogeneous complex
Monge-Ampere equations are closely related to foliations by
holomorphic curves (cf. \cite{Lempt83}, \cite{Semmes92},
\cite{dona01}). In \cite{Semmes92}, S. Semmes formulated the
Dirichlet problem for (\ref{eq:hcma1}) in terms of a foliation by
holomorphic curves with boundary in a totally real submanifold in
the complex cotangent bundle of the underlying manifold.
%Donaldson rediscovered this formulation and used it to study
%deformation of smooth solutions for (\ref{eq:hcma1}) when the
%boundary value varies. Our work was also inspired by
%\cite{semmes93} and \cite{dona01}. Our new contribution is to
%establish existence of foliations by holomorphic disks with
%relatively mild singularity and use them to obtain regularity of
%solutions outside a small subset. As in the case of harmonic maps
%etc., we call it partial regularity of the solutions. It seems to
%be the first time to use singular foliations systematically to
%study partial regularity for homogeneous complex Monge-Ampere
%equations. Furthermore, we have significant geometric applications
%as we stated above.
Let us first recall Semmes' construction. We associate a
hyperK\"ahler manifold $\cW_{[\omega]}$ to each K\"ahler class
$[\omega]$: Let $\{U_i\}$ be a covering of $M$ such that
$\omega|_{U_i} = \sqrt{-1}\partial\overline{\partial} \rho_i$, we
identify $(x,v_i)\in T^*U_i$ with $(y, v_j)\in T^*U_j$ if $x=y\in
U_i\cap U_j$ and $v_i=v_j + \partial (\rho_i-\rho_j)$, then
$\cW_{[\omega]}$ consists of all these equivalence classes of
$[x,v_i]$. There is an natural map $\pi: \cW_{[\omega]}\mapsto
T^*M$, assigning $(x,v_i)\in T^*U_i$ to $(x, v_i - \partial
\rho_i)$. Then the complex structure on $T^*M$ pulls back to a
complex structure on $\cW_{[\omega]}$ and there is also a
canonical holomorphic 2-form $\Omega$ on $\cW_{[\omega]}$, in
terms of local coordinates $z_\alpha, \xi_\alpha$ ($\alpha =
1,\cdots, n$) of $T^*U_i$,
$$\Omega = dz_\alpha \wedge d\xi_\alpha .$$
Now for any $\varphi\in \cH_{[\omega]}$, we can associate a
complex submanifold $\Lambda_\varphi$ in $\cW_{[\omega]}$: For any
open subset $U$ on which $\omega$ can be written as
$\sqrt{-1}\partial\overline
\partial \rho$, we define $\Lambda_\varphi|_U = $ to be the graph of
$\partial (\rho + \varphi)$. Clearly, this $\Lambda_\varphi$ is
independent of the choice of $U$. A straightforward computation
shows
\begin{equation}
\lab{eq:sympl-lang} \Omega|_{\Lambda_\varphi} = -\sqrt{-1}
\omega_\varphi,
\end{equation}
that is, ${\rm Re}(\Omega)|_{\Lambda_\varphi}=0$ and $-{\rm
Im}(\Omega)|_{\Lambda_\varphi} = \omega_\varphi > 0$. This means
that ${\Lambda_\varphi}$ is an exact Lagrangian symplectic
submanifold of $\cW_{[\omega]}$ with respect to $\Omega$.
Conversely, given an exact Lagrangian symplectic submanifold
$\Lambda$ of $\cW_{[\omega]}$, we can construct a smooth function
$\varphi$ such that $\Lambda=\Lambda_\varphi$. Hence, K\"ahler
metrics with K\"ahler class $[\omega]$ are in one-to-one
correspondence with exact Lagrangian symplectic subamnifolds in
$\cW_{[\omega]}$.

Let $\psi$ be a smooth function on $\partial \Sigma \times M$ such
that $\psi(\tau,\cdot)\in \cH_{[\omega]}$ for any $\tau\in
\partial\Sigma$. Define
\begin{equation}
\lab{eq:boundaryvalue} {\bf \Lambda}_\psi = \{ (\tau, v) \in
\partial \Sigma\times \cW_{[\omega]}~|~v\in {\bf \Lambda}_{\psi(\tau,
\cdot)}~\}.
\end{equation}
One can show that ${\bf \Lambda}_\psi$ is a totally real
submanifold in $\Sigma\times\cW_{[\omega]}$. Now let us recall a
result from \cite{Semmes92} and \cite{dona01}.
\begin{prop}
\lab{prop:foliation} Assume that $\Sigma$ is simply connected.
There is a solution $\phi$ of (\ref{eq:hcma1}) if and only if
there is a smooth family of holomorphic maps $h_x:\Sigma\mapsto
\cW_{[\omega]}$ parametrized by $x\in M$ satisfying: (1)
$\pi(h_x(z_0))=x$, where $z_0$ is a given point in
$\Sigma\backslash \partial \Sigma$; (2) $h_x(\tau)\in {\bf
\Lambda}_{\psi(\tau,\cdot)}$ for each $\tau \in \partial \Sigma$
and $x\in M$; (3) For each $z\in \Sigma$, the map
$\gamma_z(x)=\pi(h_x(z))$ is a diffeomorphism of $M$.
\end{prop}

For the readers' convenience, let us explain briefly its proof.
Let $\phi$ be a solution of (\ref{eq:hcma1}) on $\Sigma\times M$
such that $\phi (z,\cdot)\in \cH$ for any $z\in \Sigma$. Define
$\cD_\phi\subset T(\Sigma\times M)$ by
\begin{equation}
\label{eq:distribution} \cD_\phi = \{{\partial \over {\partial z}}
+ v\in T_{(z,p)}(\Sigma\times M)~|~i_{{\partial \over {\partial
z}} + v}(\pi_2^*\omega +\partial\overline{\partial}
\phi)=0~\},~~~(z,p)\in \Sigma\times M.
\end{equation}
Then $\cD$ is a holomorphic integrable distribution. If $\Sigma$
is simply-connected and $\phi(z,\cdot)\in \cH$ for each $z\in
\Sigma$, then the leaf of $\cD$ containing $(z_0,x)$ is the graph
of a holomorphic map $f_x:\Sigma\mapsto M$ with $f_x(z_0)=x$. If
we write $f_x(z)=\sigma_z(x)$ we get a family of diffeomorphisms
$\sigma_z$ of $M$ with $\sigma_{z_0}={\rm Id}_M$. Now for any
fixed $z$ we have a K\"ahler form $\omega +
\sqrt{-1}\partial\overline\partial \phi(z,\cdot)$ on $M$ and hence
a section $s_z:M\mapsto \cW_{[\omega]}$ whose image is an exact
Lagrangian symplectic graph $\Lambda_{\phi(z,\cdot)}$. Then
$h_x(z)=\gamma_z(x)=s_z(f_x(z))$ as required. This process can be
reversed.

In \cite{dona01}, S. Donaldson used this fact to study deformation
of smooth solutions for (\ref{eq:hcma1}) when the boundary value
varies. Theorem \ref{th:almostsmoothsolution} is proved by
establishing existence of foliations by holomorphic disks with
relatively mild singularity, more precisely, we will show that for
a generic boundary value, there is an open set in the {\it moduli}
space of holomorphic discs which generates a foliation on
$\Sigma\times M\backslash S$ for a closed subset $S$ of
codimension at least one.

Now let us fix a generic boundary value $\psi$ and study the
corresponding {\it moduli} $\cM_\psi$ of holomorphic discs. First
it follows from the Index theorem that the expected dimension of
this {\it moduli} is $2n$. Recall that a holomorphic disc $u$ is
regular if the linearized $\overline\partial$-operator
$\overline\partial _u$ has vanishing cokernel. The {\it moduli}
space is smooth near a regular holomorphic disc. Following
\cite{dona01}, we call $u$ super-regular if there is a basis
$s_1,\cdots ,s_{2n}$ of the kernel of $\overline\partial _u$ such
that $d\pi(s_1)(x),\cdots, d\pi(s_{2n})(x)$ span $T_{u(x)}M$ for
every $x\in \Sigma$, where $\pi: \cW_{[\omega]} \mapsto M$ is the
natural projection. We call $u$ almost super-regular if
$d\pi(s_1)(x),\cdots, d\pi(s_{2n})(x)$ span $T_{u(x)}M$ for every
$x\in \Sigma\backslash \partial \Sigma$. Clearly, the set of
super-regular discs is open.

One of our crucial observations is that Semmes' arguments can be
made local along super-regular holomorphic discs.
\begin{theo}
\lab{th:polarized=almostsuper} For a generic boundary value
$\psi$, an almost smooth solution of (\ref{eq:hcma1}) corresponds
to a nearly smooth foliation which can be described as follows:
There is an open subset $\cU_\psi\subset\cM_\psi$ of super-regular
discs such that the images of these discs in $\Sigma\times M$ give
rise to a foliation on an open-dense set of $\Sigma\times M$
 such that
\begin{enumerate}
\item this foliation can be extended to be a continuous foliation
by holomorphic discs in an open set $\tilde\cV_{\phi_0}\subset
\Sigma_0\times M $ such that it admits a continuous lifting in
$\Sigma\times \cW_M$;
%satisfying:
%\begin{enumerate} \item it The complement of each $\cV_{\phi_0}\cap
%(\{z\}\times M)$ in $M$ has codimension at least one; \item
\item the complement of $\tilde \cV_{\phi_0}$ in $\Sigma_0\times
M$ is locally extendable. \item The leaf vector (cf. equation
\ref{eq:distribution}) induced by the foliation in $\cV_{\phi_0}$
is uniformly bounded.
\end{enumerate}
%Moreover, $\bar U$ is a smooth closed subset of $\cM_\psi$ and the evaluation map from $\Sigma
%\times \bar U$ to $\Sigma\times M$ is a continuous onto map with degree 1.
\end{theo}

We can prove that these locally constructed solution are actually
compatible to each other if their domains overlap. If the set of
super regular discs is $M$, then these locally constructed
solutions give rise to a smooth solution to (\ref{eq:hcma1}). In
general, one gets only a solution to (\ref{eq:hcma1}) in
$(\Sigma\times M) \setminus S$ with appropriate boundary condition
on $(\partial \Sigma \times M)\setminus S$. We can apply the
Maximum Principle along super-regular leaves to get a uniform
$C^{1,1}$-bound $\phi$ on $(\Sigma \times M)\setminus S$. This
uniform $C^{1,1}$-bound can be used to get a solution by patching
together local solutions to (\ref{eq:hcma1}) along super-regular
leaves.

Theorem \ref{th:almostsmoothsolution}  will follow from the
following

\begin{theo}
\lab{th:existenceofpartialfoliation} For a generic boundary value
$\psi$, there is a nearly smooth foliation generated by an open
set of the corresponding {\it moduli} space $\cM_\psi$. Definition
of {\it nearly smooth foliation} is already given in the statement
of Theorem 1.7.
% Moreover, the singular set of this nearly smooth foliation is of codimension at least one.
\end{theo}

Now we outline the proof of Theorem
\ref{th:existenceofpartialfoliation}. Let $\psi$ be a generic
boundary value such that $\cM_\psi$ is smooth. This follows from a
result of Oh  \cite{Oh952} on transversality. By a similar (but
different) transversality argument, one can show that there is a
generic path $\psi_t$ ($0\le t\le 1$) such that $\psi_0=0$ and
$\psi_1=\psi$ and the total {\it moduli} space $\tilde \cM =
\bigcup_{t\in [0,1]} \cM_{\psi_t}$ is smooth, moreover, we may
assume that $\cM_{\psi_t}$ are smooth for all $t$ except finitely
many $t_1,\cdots, t_N$ where the moduli space may have isolated
singularities. It follows from Semmes and Donaldson
(\ref{prop:foliation}) that $\cM_0$ has at least one component
which gives a foliation for $\Sigma\times M$. We want to show that
this component will deform to a component of $\cM_\psi$ which
generates a nearly smooth foliation. We will use the continuity
method. Assume that $\phi$ is the unique $C^{1,1}$-solution of
(\ref{eq:hcma1}) with boundary value $\psi_t$ for some $t\in
[0,1]$. Let $f$ be any holomorphic disc in the component of
$\cM_{\psi_t}$ which generates the corresponding foliation.

\begin{lemm}
\lab{le:areabound} There is a uniform upper bound on area of $f$,
where the area is with respect to the induced metric on
$\cW_{[\omega]}$ by $dzd\bar z $ on $\Sigma$ and $\omega$ on $M$.
\end{lemm}
\begin{proof} First we observe
\begin{equation}
\lab{eq:areabound} {\rm Area} f(\Sigma)\le C {\sqrt{-1}}
\int_\Sigma f^*(dz\wedge d\bar z + \omega),
\end{equation}
where $C$ is a constant depending only on $C^{1,1}$-norm of
$\phi$. By direct computations using (\ref{eq:hcma1}), we have
$\partial \overline{\partial }\phi(z,f(z)) = - f^*\omega$. Then,
integrating by parts, we can bound the area of $f$ in terms of the
area of $\Sigma$ and the $C^{1,1}$-bound of $\phi$.
\end{proof}
By Gromov's compactness theorem, any sequence of holomorphic discs
with uniformly bounded area has a subsequence which converges to a
holomorphic disc together with finitely many bubbles. These
bubbles which occur in the interior are holomorphic spheres, while
bubbles in the boundary might be holomorphic spheres or discs. We
will show that no bubbles can actually occur.

For a fixed totally real submanifold, holomorphic bubbles can not
occur in the boundary since boundaries of discs lie in a fixed
totally real submanifold. If a sequence of totally real
submanifolds converges to a given totally real submanifold, there
are two limiting processes, one concerns how fast the bubbles form
and move to the boundaries of discs, while the other is about how
fast the sequence of totally real submanifolds approaches to the
limiting submanifold. The uniform $C^{1,1}$ bound on $\phi$ can be
used to show that the two limiting processes are exchangeable.
Consequently, one can show that there are no bubbles along
boundary.

We can also rule out bubbles in the interior of discs.
Heuristically speaking, an interior bubble corresponds to a
holomorphic map from $S^2$ into the target manifold $\cH$.
According to E. Calabi and X. Chen \cite{chen002}, this infinite
dimensional space $\cH$ is non-positively curved in the sense of
Alexanderov, consequently, there are no non-trivial holomorphic
spheres in $\cH_{[\omega]}$! This heuristical argument implies
that there are no interior bubbles. Indeed, there is a rigorous
proof for this fact. The proof is much more involved and will be
presented elsewhere.

Since there are no bubbles arised either in the boundary or
interior of the disc $\Sigma$, the Fredholm index of holomorphic
discs is invariant in limiting process. This is an important fact
needed in our doing deformation theory.

In order to get a nearly smooth foliation, we need to prove that
the {\it moduli} space has an open set of super-regular
holomorphic discs for each $t$. First we observe that the set of
super-regular discs is open. Moreover, using the transversality
arguments, one can show that for a generic path $\psi_t$, the
closure of all super-regular discs in each $\cM_{\psi_t}$ is
either empty or forms an irreducible component. It implies the
openness. It remains to proving that each moduli has at least one
super-regular disc. It is done by using capacity estimate which we
explain briefly in the following.

%We have seen in the above
%\begin{equation}
%\lab{eq:equationonphi}\partial_z\overline{\partial_z } \phi(z,f(z)) = - f^*\omega.
%\end{equation}

%As discussed before, a smooth solution to HCMA equation gives rise to a foliation of holomorphic discs
%in $\Sigma\times M.\;$
Consider the bundle $\cE=\pi_2^*TM$ over $\Sigma\times M$. Each
almost smooth solution $\phi$ of (\ref{eq:hcma1}) induces an
Hermitian metric on $\cE|_{{\cal R}_\phi}$, where ${\cal R}_\phi$
was defined in Theorem \ref{th:almostsmoothsolution}. If $f$ is a
super-regular disc, then $\cE$ pulls back to an Hermitian bundle
over $\Sigma$ with fiber $T_{f(z)}M$ and metric
$\omega_{\phi(z,\cdot)}(f(z))$ over $z\in \Sigma$. It turns out
that the curvature of this Hermitian bundle is non-positive. This
fact plays a crucial role in our work. More precisely, we have

\begin{lemm}
\lab{lemm:curvatureformula} Let $\phi$ be a solution of
(\ref{eq:hcma1}) and $f$ be a super-regular holomorphic disc as
above, then the curvature form $F$ of the metric $g_\phi$
described above is given by
\[
g_\phi(F(u), v) = -g_\phi(u({\overline{\partial_z f}}),
{v}(\overline{\partial_z f}) ) ~~~~u, v \in TM.
\]
In particular, the curvature is non-positive. Moreover, the
foliation is holomorphic along $f$ if and only if the curvature
vanishes.
\end{lemm}

The determinant $\wedge^n \cE$ restricts to an Hermitian line
bundle over any given super-regular disc. The corresponding
Hermitian metric, denoted by $f^*\omega^n_\phi$, at $z$ is
$\omega^n_{\phi(z,\cdot)}(f(z))$. An immediate corollary of above
lemma is that the curvature of this line bundle is non-positive
\footnote{This fact was first proved by S. K. Donaldson and X. X.
Chen independently (via different methods) while both of them were
visiting at Stanford University.}. Moreover, there are constants
$C_1, C_2 $ which depend only on the background metric $\omega$
such that
\begin{equation}
\label{eq:volumeratio1} \Delta \left(\log {{f^*\omega_\phi^n}\over
f^*\omega^n} + C_1 \varphi\right) \geq 0,
\end{equation}
and
\begin{equation}
 \label{eq:volumeratio2}
\Delta \left(\log {{f^*\omega_\varphi^n}\over f^*\omega^n} + C_2
\varphi\right) \leq - {\rm tr}(F),
\end{equation}
where $\Delta$ denotes the standard Laplacian operator on $\Sigma$
and $\varphi(z)=\phi(z,f(z))$. It follows that $ \log
{{f^*\omega_\phi^n}\over f^*\omega^n} + C_1 \varphi$ is
subharmonic and uniformly bounded on the boundary $\partial
\Sigma$. The $C^{1,1}$-estimate in \cite{chen991} implies that
this function is uniformly bounded from above. Moreover, the
difference of two functions $\log {{f^*\omega_\phi^n}\over
f^*\omega^n} + C_1 \varphi$ and $\log {{f^*\omega_\phi^n}\over
f^*\omega^n} + C_2 \varphi$ is uniformly bounded. In addition, we
have \begin{equation}- \Delta{\rm tr}(F)\ge \frac{2}{n} (-{\rm
tr}(F))^2 \label{eq:thirdderivatives}\end{equation} for some
positive constant $c$. Following \cite{Osserman57} and
\cite{Calabi57}, one can use this differential inequality to
derive an interior estimate on ${\rm tr}(F)$ (details will appear
elsewhere). Applying this estimate on ${\rm tr}(F)$ to the above
equations, we can derive an Harnack-type inequality
${f^*\omega^n_\phi \over f^*\omega^n}$ in the interior of
$\Sigma$.

Now let us introduce the notion of Capacity for super-regular
holomorphic discs:
\begin{defi} For any super-regular disc $f$ in an moduli space ${\cal M}_{\psi}$, we define
its capacity by
\[
{Cap}(f) = {\sqrt{-1}\over 2} \;\displaystyle\int_\Sigma\;
{f^*\omega^n\over f^*\omega_\phi^n} \; dz\wedge d\bar z.
\]
\end{defi}

Using the Harnack-type inequality mentioned above, one can control
the lower bound of ${f^*\omega_\phi^n\over f^*\omega^n}$ in the
interior of $\sigma$ in terms of upper bound of the capacity of
$f$. This has an important corollary for compactness of
super-regular discs with uniformly bounded capacity.

\begin{theo} Let $f_i$ be any sequence of super-regular discs
in $\cM_{\psi_{t_i}}$ which converges smoothly to an embedded disc
$f_\infty$ in $\cM_{\psi_{t_\infty}}$. If the capacities
$Cap(f_i)$ are uniformly bounded, then the limiting disc
$f_\infty$ is also super-regular.
\end{theo}

In fact, Lemma \ref{lemm:curvatureformula} was already needed when
we extended Semmes' correspondence to almost smooth solutions of
(\ref{eq:hcma1}) and nearly smooth foliations. For this local
extension, we first construct smooth solutions of (\ref{eq:hcma1})
along super-regular leaves and then glue them together to a
solution $\phi$ on an open and dense subset $V_{good}\subset
\Sigma\times M$, but we need to establish $C^{1,1}$-bound of
$\phi$. Once this bound is established, the maximum principle
implies that $\phi$ coincides with the solution in \cite{chen991}.
The $C^{1,1}$-bound of $\phi$ follows from the following

\begin{theo}
\lab{th:sectionbound} For any global holomorphic section
$s:\Sigma\rightarrow \cE$, the norm of $s$ with respect to
$g_\phi$ achieves its maximum value at the boundary of the disc.
\end{theo}
%\bibliography{test}

\begin{thebibliography}{10}

\bibitem{Au76}
T.~Aubin,
\newblock Equations du type de Monge-Ampere sur les varietes K\"ahleriennes compactes.
\newblock {\em C. R. Acad. Sci. Paris}, 283 (1076), 119--121.

\bibitem{Bando87}
S.~Bando and T.~Mabuchi,
\newblock {U}niquness of {E}instein {K}\"ahler metrics modulo connected group actions.
\newblock In {\em Algebraic {G}eometry}, Advanced Studies in Pure Math., 1987, 11--40.

\bibitem{Bedford76}
E.D. Bedford and T.A. Taylor,
\newblock The {D}richelet problem for the complex {M}onge-{A}mpere operator.
\newblock {\em Invent. Math.}, 37 (1976), 1--44.


\bibitem{Calabi57}
E.~Calabi,
\newblock An extension of E. Hopf's maximum principle with an application to Riemannian geometry.
\newblock {\em Duke. J. Math.}, 25 (1957), 45-56.

\bibitem{calabi82}
E.~Calabi,
\newblock Extremal {K}\"ahler metrics.
\newblock In {\em Seminar on Differential Geometry}, volume~16 of {\em 102},
259--290. Ann. of Math. Studies, University Press, 1982.

\bibitem{chen991}
X.~X. Chen,
\newblock Space of {K}\"ahler metrics.
\newblock {\em J. Diff. Geom.}, 56 (2000), 189--234.

\bibitem{chen002}
E.~Calabi and X.~X. Chen,
\newblock Space of {K}\"ahler metrics. II.
\newblock {\em J. Diff. Geom.}, 61 (2002), 173--193.


\bibitem{Dona96}
S.K. Donaldson,
\newblock Symmetric spaces, K\"ahler geometry and {H}amiltonian dynamics.
\newblock {\em Amer. Math. Soc. Transl}, Ser. 2, 196 {1999), 13--33.

\bibitem{dona01}
S.K. Donaldson,
\newblock Holomorphic {D}iscs and the complex {M}onge-{A}mpere equation, 2001.
\newblock {\em {J}ournal of {S}ymplectic {G}eometry}, 1 (2000), 171-196.

\bibitem{dona02}
S.K. Donaldson,
\newblock Scalar curvature and projective embeddings. I.,
\newblock {\em J. Diff. Geom.}, 59 (2001), 479--522.


\bibitem{Lempt83}
L.~Lempert,
\newblock Solving the degernerate {M}onge-{A}mpere equation with one
  concentrated singularity.
\newblock {\em Math. Ann.}, 263 (1983), 515--532.

\bibitem{Ma87}
T.~Mabuchi,
\newblock Some {S}ymplet}ic geometry on compact k\"ahler manifolds {I}.
\newblock {\em Osaka {J}. {M}ath.}, 24 (1987), 227--252.

\bibitem{Oh952}
Yong-Geun Oh.
\newblock Riemann-{H}ilbert problem and application to the perturbation theory
  of analytic discs.
\newblock {\em KYUNGPOOK Math., J.}, 35:38--75, 1995.


\bibitem{Osserman57}
R. ~Osserman,
\newblock On the inequality $\triangle u\geq f(u)$.
\newblock {\em Pacific J. Math.}, 7 (1957), 1641-1647.

\bibitem{Semmes92}
S.~Semmes,
\newblock Complex Monge-Ampere and sympletic manifolds.
\newblock {\em Amer. J. Math.}, 114 (1992), 495--550.

%%\bibitem{semmes93} R.R. Coifman and S.~Semmes.
%%\newblock Interpolation of {B}anach spaces, {P}erron process, and {Y}ang {M}ills.
%%\newblock {\em Amer. J. Math.}, 115(2):243--278, 1993.

\bibitem{Tian97}
G.~Tian,
\newblock K\"ahler-{E}instein metrics with positive scalar curvature.
\newblock {\em Invent. Math.}, 130 (1997), 1--39.

\bibitem{tian001}
G.~Tian,
\newblock Bott-Chern forms and geometric stability
\newblock {\em Discrete Contin. Dynam. Systems}, 6 (2000), 1--39.

\bibitem{paultian04}
S. Paul and G. Tian,
\newblock Algebraic and Analytic K-Stability, preprint, 2004.

\bibitem{tianzhu00}
G.~Tian and X.~H. Zhu,
\newblock Uniqueness of KŠhler-Ricci solitons.
\newblock {\em Acta Math.}, 184 (2000), 271--305.


\bibitem{Yau78}
S.~T. Yau,
\newblock On the {R}icci curvature of a compact {K}\"ahler manifold and the
  complex {M}onge-{A}mpere equation, ${I}^*$.
\newblock {\em Comm. Pure Appl. Math.,}, 31 (1978), 339--441.

\end{thebibliography}

\end{document}